\documentclass[12pt]{article}
\usepackage{amsmath,amsfonts,amssymb,amsthm,graphicx,float,color,caption,verbatim}
\usepackage{mathrsfs}\usepackage{hyperref}

\usepackage{geometry}
\theoremstyle{plain}
\newtheorem{theorem}{Theorem}[section]

\theoremstyle{definition}
\newtheorem{definition}[theorem]{Definition}

\theoremstyle{remark}

\title{Inverse Problems of a Fractional Differential Equation with Bessel Operator}
\author{Fatma Al-Musalhi, Nasser Al-Salti, Sebti Kerbal}
\begin{document}
\maketitle
\begin{abstract}
Inverse initial and inverse source problems of a time-fractional differential equation with Bessel operator are considered. Results on existence and uniqueness of solutions to these problems are presented. The solution method is based on series expansions using a set of Bessel functions of order zero. Convergence of the obtained series solutions is also discussed. 
\end{abstract}
\section{Introduction}
Inverse problems have been a very active field of research in mathematics and science. They appear in modeling a wide variety of real life problems, e.g., imaging, Magnetic Resonance Imaging, Computerized Tomography, signal processing and many other applications (see, for example,\cite{BB}, \cite{uhlmann1} and \cite{uhlmann2}). Inverse problems come into various types, for example, inverse initial problems where initial data are unknown and inverse source problems where the source term is unknown. These unknown terms are to be determined using extra boundary data.

Fractional differential equations, on the other hand, become an important tool in modeling many real-life problems and hence there has been growing interest in studying inverse problems of time fractional equations \cite{ams}, \cite{furati}, \cite{kirane and malik}, \cite{zhang}. For example, in \cite{zhang}, Zhang and Xu studied inverse source problem for a fractional diffusion equation where solutions are found based on the method of eigen-function expansion. 

In \cite{initial}, authors considered initial inverse problem in heat equation with Bessel operator. They expressed the solution of the problem and the initial temperature distribution in terms of orthogonal set of Bessel functions. These types of functions arises in modeling of chemical engineering process including hydrodynamics, bio-processes, diffusion and heat transfer (see, for example, \cite{initial},\cite{Pet}). Properties of Bessel functions are given in \cite{tolstov}, \cite{watson}.

This paper concerns with the existence and uniqueness of solutions to inverse initial and inverse source problems for a time fractional differential equation with Bessel operator. Solutions are obtained in the form of series expansion using an orthogonal set of Bessel functions. 

\section{Preliminaries }
In this section, we recall some definitions and results which will be used later in this article.\\
First, Caputo fractional derivative  of order $0<\alpha<1$, $^cD_{0|t}^{\alpha}f(t)$, is defined by (\cite{gm}, p. 228):

$$ ^cD_{0|t}^{\alpha}f(t)= \dfrac{1}{\Gamma(1-\alpha)}\displaystyle\int_{0}^{t}\frac{f' (\tau)}{(\tau-t)^{\alpha}}d\tau,
$$
where $\Gamma(.)$ denotes the well-known gamma function.
The  Laplace transform of Caputo fractional derivative of order $0<\alpha<1$ is
$$
\mathscr{L}\lbrace ^cD_{0|t}^{\alpha}f(t)\rbrace (s)=s^\alpha F(s)-s^{\alpha-1} f(0^+),\\
$$ where $F(s)=\mathscr{L}\left\lbrace f(t);s \right\rbrace.$ 
Also, we need Mittag-Leffler function which plays a major role in describing solutions to fractional order differential equations. The Mittag-Leffler function with one parameter, $E_{\alpha}(z)$, is defined as
$$
E_{\alpha}(z)=\sum_{k=0}^{\infty}\dfrac{z^k}{\Gamma(\alpha k +1)},\,\ \text{Re}(\alpha)>0,\, z\in C.
$$
If $z=\lambda t^{\alpha}$, the Laplace transform of $E_{\alpha}(\lambda t^\alpha)$ is given by
$$
\mathscr{L
}\left\lbrace E_{\alpha}(\lambda t^\alpha)\right\rbrace=\dfrac{s^{\alpha-1}}{s^{\alpha}-\lambda}.
$$ 
Moreover, Mittag-Leffler function $E_{\alpha}(\lambda t^\alpha)$ is bounded, i.e, (see \cite{Prabhakar})
\begin{equation}\label{bm}
E_{\alpha}(\lambda t^\alpha)\leq M.
\end{equation}
Here and throughout this paper, $M$ denotes a positive constant. 
Now, we  introduce Bessel's equation of order $\nu\geq 0$ :
\begin{equation}
y''+\dfrac{1}{x}y'+ \left(\lambda^2-\dfrac{\nu^2}{x^2}\right)y=0,
\end{equation}
which when $\nu$ is integer has the following solution:
$$y(x)=C_1 J_{\nu}(\lambda x)+C_2 Y_{\nu}(\lambda x),$$
where $J_{\nu}(x)$ and $Y_{\nu}(x)$ are Bessel functions of the first and second kind of order $\nu$, respectively. The Bessel functions of the second kind are not bounded near the point $x=0$ and hence for a bounded solution near zero, we need to have $C_2=0$, i.e., 
$$ y(x)= C_1 J_{\nu}(\lambda x).$$
Moreover, if the boundary condition $y(1)=0$ is imposed, then the parameter $\lambda$ must satisfy $J_{\nu}(\lambda)=0$, i.e., the values of $\lambda$ represent the zeros of the Bessel function $J_{\nu}(x)$, which has the following asymptotic representation (see \cite{tolstov}, page 213)
\begin{equation*}
J_{\nu}(x)=\sqrt{\frac{2}{\pi x}}\sin\left(x-\frac{\nu\pi}{2}+\frac{\pi}{4}\right)+\frac{r_{\nu}(x)}{x\sqrt{x}},
\end{equation*}
where the function $r_{\nu}(x)$ remains bounded as $x\rightarrow \infty$. Hence, for arbitrary large $k$, the zeros of  $J_{\nu}(x) $ are given by (see, for example, \cite{tolstov}, page 214)
\begin{equation*}
k\pi+\frac{\nu\pi}{2}-\frac{\pi}{4}.
\end{equation*}
  
Next, we define a Fourier-Bessel expansion of a given function $f(x)$ as follows:
\begin{definition}( see for example \cite{tolstov}, page 211 and \cite{watson}, sec. 18.11) Let $\lambda_k, k=1,2,\cdots$, be the positive zeros of $J_{\nu}(x) $ arranged in ascending order of magnitude. Then, the expansion of a function $f(x)$ in terms of the Bessel functions $J_{\nu}(\lambda_k x), k=1,2,\cdots$, known as a Fourier-Bessel series, is given by
\begin{equation}
f(x)=\sum_{k=1}^{\infty} c_k J_\nu (\lambda_k x),
\end{equation}
where $\nu\geq -\frac{1}{2}$ and $c_k$, which are called the Fourier-Bessel coefficients for $f(x)$, are given by
\begin{equation}
c_k=\dfrac{2}{J^2_{\nu+1} (\lambda_{k})}\int_{0}^{1} x f(x) J_{\nu}(\lambda_{k} x) dx, \quad  k=1,2,\cdots.
\end{equation}
\end{definition}

Finally, some results on the uniform convergence of Fourier-Bessel series are stated in the following theorems:

\begin{theorem}( for proof, see \cite{tolstov}, pp. $230 - 231$) \label{ck}
Let $f(x)$ be a function defined on the interval $[0,1]$ such that $f(x)$ is differentiable $2s$ times ($s\geq 1$) and  
\begin{enumerate}
\item $f(0)=f'(0)=\cdots=f^{(2s-1)}(0)=0$;
\item $f^{(2s)}(x)$ is bounded ( this derivative may not exist at certain points);
\item $f(1)=f'(1)=\cdots=f^{(2s-2)}(1)=0.$\end{enumerate}
Then, the following inequality is satisfied by the Fourier-Bessel coefficient of $f(x)$:
\begin{equation}
|c_{k}|\leq \dfrac{M}{\lambda_{k}^{2s-(1/2)}}. 
\end{equation} 
\end{theorem} 

\begin{theorem}(for proof, see \cite{tolstov}, p. $225$) \label{c}
If $\nu\geq 0$ and if 
\begin{equation}
|c_k|\leq \dfrac{M}{\lambda^{1+\varepsilon}},
\end{equation}
where $\varepsilon$  is a positive constant, then the series 
$$ \sum_{k=1}^{\infty}c_k J_\nu(\lambda_k x) ,$$
converges absolutely and uniformly on $[0,1].$
\end{theorem}
For more theorems on the convergence of Fourier-Bessel series, we refer the reader to \cite{tolstov} and \cite{watson}.

The remaining two sections are the main sections of this paper, where we present the existence and uniqueness of solutions to inverse initial and inverse source problems for a time fractional differential equations with Bessel operator, respectively. The main result for each inverse problem is summarized at the end of each section.

\section{Inverse initial problem}
Here, we consider the following inverse initial problem:\\ Find a pair of functions $\left\lbrace u(x,t), g(x) \right\rbrace $, in a  rectangular domain $$\Omega=\left\lbrace (x,t):0<x<1,\; 0<t<T\right\rbrace,$$ which satisfies the fractional differential equation
\begin{equation}\label{prbu}
^c D^{\alpha}_{0|t}u=\dfrac{\partial^2 u}{\partial x^2}+\dfrac{1}{x}\dfrac{\partial u}{\partial x},
\end{equation}
along with the boundary conditions 
\begin{equation} \label{BC1}
\underset{x\rightarrow 0}{\lim} \;x \;\dfrac{\partial u}{\partial x}=0,\quad u(1,t)=0,
\end{equation}
the initial condition
\begin{equation}\label{ICg}
u(x,0)=g(x),
\end{equation}
and the overdetermining condition
\begin{equation}\label{condf}
u(x,T)=f(x).
\end{equation}
Here $f(x)$ is a given function and  $ \,^cD^{\alpha}_{0|t}$ stands for the Caputo fractional derivative of order $0<\alpha< 1$.
\subsection{Solution Method}
Using separation of variables method, we obtain, from (\ref{prbu}) and (\ref{BC1}), the following  spectral problem
\begin{eqnarray}\label{sepr} 
&X''+\dfrac{1}{x}X'+\lambda^2 X=0,& 0<x<1, \\
\label{besselcond}
&\underset{x\rightarrow0}{\lim}\; x X'(x)=0,&X(1)=0, 
\end{eqnarray}
which is a self-adjoint problem. Note that the above spectral equation $(\ref{sepr})$ represents Bessel's equation of order zero and hence one can show that the eigenfunctions, $X_k$, for the problem $(\ref{sepr})$ - $(\ref{besselcond})$ are the Bessel functions of the first kind of order zero, i.e.,
$$X_k=J_0(\lambda_k x),\quad k=1,2,3,\cdots.$$
The corresponding eigenvalues, $\lambda_k$, are the positive zero of $J_0(x)$, i.e., they satisfy $J_0(\lambda_k)=0$ and for large $k$, they are given by 
$$\lambda_k = k\pi- \dfrac{\pi}{4}$$
as discussed earlier. The set of eigenfunctions $\{X_k \}$ forms a complete orthogonal basis in $L^{2}(0,1)$(\cite{hig}, p. 40). 
Hence, the functions $u(x,t)$, $g(x)$ and  $f(x)$ can be expanded as follows:
\begin{equation}\label{solu}
u(x,t)=\sum_{k=1}^{\infty} u_k(t) J_0(\lambda_k x),
\end{equation}
\begin{equation} \label{solg}
g(x)=\sum_{k=1}^{\infty} g_k J_0(\lambda_k x), \quad f(x)=\sum_{k=1}^{\infty} f_k J_0(\lambda_k x),
\end{equation}
where  $u_k(t),$ $g_k$ are unknowns and the coefficient $f_k$ is given by 
\begin{equation}
f_k=\dfrac{2}{J_{1}^{2}(\lambda_k)}\int_{0}^{1}x\, f(x) J_{0}(\lambda_kx) dx.
\end{equation}
Substituting $(\ref{solu})$ and $(\ref{solg})$ into $(\ref{prbu})$, $(\ref{ICg})$ and $(\ref{condf})$ gives the following linear  fractional differential equation: 
\begin{equation}\label{fde}
 \,^c D^{\alpha}_{0|t}u_k(t)+\lambda_k^{2}  u_k(t)=0,\end{equation}
along with the conditions: 
$$u_k(0)=g_k,\quad u_k(T)=f_k .$$
Applying Laplace transform to equation $(\ref{fde})$, we obtain (\cite{gm}, p. 243)
$$U_k(s)=\dfrac{s^{\alpha-1}}{s^{\alpha}+\lambda_k^2} g_k,$$
and consequently, we get 
$$u_k(t)= g_k E_{\alpha}(-\lambda_k^2 t^\alpha).$$
Next, using the  condition $u_k(T)=f_k $, we get  
$$g_k=\dfrac{f_k}{E_{\alpha}(-\lambda_k^2 T^\alpha)}.$$
Therefore, $u(x,t)$ and $g(x)$ can now be written as
\begin{equation}\label{seriesu}
 u(x,t) =\sum_{k=1}^{\infty}   \dfrac{f_k}{E_{\alpha}(-\lambda_k^2 T^\alpha)} E_{\alpha}(-\lambda_k^2 t^\alpha)J_0(\lambda_k x),\end{equation}
\begin{equation}g(x)=\sum_{k=1}^{\infty} \dfrac{f_k}{E_{\alpha}(-\lambda_k^2 T^\alpha)} J_0(\lambda_k x).
\end{equation}
To complete the existence of solution, we need to show that the series representations of $u(x,t),$ $g(x)$, $\,^c D^{\alpha}_{0|t}u(x,t),$ and $u_{xx}(x,t)$ converge uniformly.
We start by estimating the coefficient of the Fourier-Bessel series of $u(x,t)$. Since $E_{\alpha}(-\lambda t^\alpha)$ is bounded, we have 
$$ \left|\dfrac{ f_k E_{\alpha}(-\lambda_k^2 t^\alpha)}{E_{\alpha}(-\lambda_k^2 T^\alpha)}\right|\leq   M \left|f_k\right|,$$
whereupon using Theorem $\ref{ck}$ for $s=1$, we obtain  
$$\left|\dfrac{ f_k E_{\alpha}(-\lambda_k^2 t^\alpha)}{E_{\alpha}(-\lambda_k^2 T^\alpha)}\right|\leq \dfrac{M}{\lambda_{k}^{3/2}}.$$
Therefore, by Theorem $\ref{c}$, we conclude that the series representation of $u(x,t)$ converges absolutely and uniformly in $\Omega$.
Similarly, one can show that the series $$g(x)=\sum_{k=1}^{\infty} \dfrac{f_k}{E_{\alpha}(-\lambda_k^2 T^\alpha)} J_0(\lambda_k x)$$ converges absolutely and uniformly in $\Omega$.\\
Now, the series expansion of $^c D^{\alpha}_{0|t} u(x,t)$ is given by 
$$\,^c D^{\alpha}_{0|t} u(x,t)=\sum_{k=1}^{\infty} \dfrac{-\lambda_k^2 f_k }{E_{\alpha}(-\lambda_k^2 T^\alpha)}E_{\alpha}(-\lambda_k^2 t^\alpha)J_0(\lambda_k x).$$
Then, for convergence we have the following estimate
$$\left|\dfrac{\lambda_k^2 f_k }{E_{\alpha}(-\lambda_k^2 T^\alpha)}E_{\alpha}(-\lambda_k^2 t^\alpha)\right|\leq M \lambda_k^2 \left|f_k\right|,$$
and by using theorem $\ref{ck}$ for $s=2$, we get 
$$M \lambda_k^2 \left|f_k\right| \leq  \dfrac{M \,\lambda_k^2}{\lambda_k^{7/2 }}= \dfrac{M}{\lambda_k^{3/2}}.$$
Therefore, the uniform and absolute convergence of the series expansion of $\,^c D^{\alpha}_{0|t} u(x,t)$ 
is ensured by Theorem $\ref{c}$. Finally, it remains to show the uniform  convergence of the series representation of $u_{xx}(x,t)$, which is given by  
$$u_{xx}(x,t)=\sum_{k=1}^{\infty} \dfrac{ f_k }{E_{\alpha}(-\lambda_k^2 T^\alpha)}E_{\alpha}(-\lambda_k^2 t^\alpha)\dfrac{d^2 }{d x^2}\left( J_0(\lambda_k x)\right).$$
Using properties of Bessel functions, namely, (see \cite{tolstov}) 
\[\begin{array}{rl}
J'_{0}( x)&=-J_{1}( x),\\
2J'_{1}(x)&=J_{0}(x) - J_{2}( x),\end{array}\]
the above representation can be rewritten as
$$
u_{xx}(x,t)=\displaystyle\sum_{k=1}^{\infty} \dfrac{ \lambda_{k}^{2} f_k }{2 E_{\alpha}(-\lambda_k^2 T^\alpha)}E_{\alpha}(-\lambda_k^2 t^\alpha)\left[J_2(\lambda_k x)- J_0(\lambda_k x)\right],$$
and hence, for convergence, we have the following estimate:
$$
\left| \lambda_{k}^{2} f_k  \dfrac{E_{\alpha}(-\lambda_k^2 t^\alpha)}{2 E_{\alpha}(-\lambda_k^2 T^\alpha)} \right| \leq M \lambda_k^2 \left|f_k\right|\leq \dfrac{M \lambda_k^2 }{ \lambda_k^{7/2}}= \dfrac{M}{ \lambda_k^{3/2}}.$$
Therefore, by Theorem $\ref{c}$, the series expansion of $u_{xx}(x,t)$ is uniformly and absolutely convergent in $\Omega$.
Note that the uniform convergence of $u_x(x,t)$  follows from the convergence  $u_{xx}(x,t).$

\subsection{Uniqueness of solution}
Suppose $\left\lbrace u_1(x,t),g_1(x) \right\rbrace $ and  $\left\lbrace u_2(x,t),g_2(x) \right\rbrace $ are two solution sets of the inverse problem (\ref{prbu}) - (\ref{condf}), then $\widehat{u}(x,t)= u_1(x,t)-  u_2(x,t)$ and $\widehat{g}(x)= g_1(x)-g_2(x)$ satisfy the following boundary value problem:
\begin{eqnarray}
^c D^{\alpha}_{0|t}\widehat{u}=\dfrac{\partial^2 \widehat{u}}{\partial x^2}+\dfrac{1}{x}\dfrac{\partial\widehat{u}}{\partial x},&(x,t)\in \Omega,\\
\underset{x\rightarrow 0}{\lim} \;x \;\dfrac{\partial \widehat{u}}{\partial x}=0,\quad \widehat{u}(1,t)=0,& 0<t<T,\\
\widehat{u}(x,0)=\widehat{g}(x),\quad \widehat{u}(x,T)=0& 0<x<1.
\end{eqnarray}
Using the completeness property of system $(\ref{sepr})$, one can then show that 
$\widehat{u}(x,t)=0$ and $\widehat{g}(x)=0,$ in $\Omega$, which in turn implies the uniqueness of solution of the inverse problem (\ref{prbu}) - (\ref{condf}).\\
The main result of this section on existence and uniqueness of solution can be summarized in the following theorem:
\begin{theorem}
Assume that the function $f$ is differentiable $4$ times such that :
\begin{itemize}
\item $f(0)=f'(0)=f''(0)=f'''(0)=0$;
\item $f(1)=f'(1)=f''(1)=0$;
\item $f^{(4)}(x)$ is bounded;\end{itemize} 
then the inverse problem 
\[\begin{array}{ll}
^c D^{\alpha}_{0|t}u=\dfrac{\partial^2 u}{\partial x^2}+\dfrac{1}{x}\dfrac{\partial u}{\partial x},& 0<x<1,\; 0<t<T,\\
\underset{x\rightarrow 0}{\lim} \;x \;\dfrac{\partial u}{\partial x}=0,\quad u(1,t)=0,&0<t<T,\\
u(x,0)=g(x),\quad u(x,T)=f(x),&0<x<1,
\end{array} \]
 has a unique solution given by
 \[\begin{array}{rl}
   u(x,t)=&\displaystyle\sum_{k=1}^{\infty}   \dfrac{f_k}{E_{\alpha}(-\lambda_k^2 T^\alpha)} E_{\alpha}(-\lambda_k^2 t^\alpha)J_0(\lambda_k x),\\
 g(x)=&\displaystyle\sum_{k=1}^{\infty} \dfrac{f_k}{E_{\alpha}(-\lambda_k^2 T^\alpha)} J_0(\lambda_k x),
\end{array}\]
where,
\[f_k=\dfrac{2}{J_{1}^{2}(\lambda_k)}\int_{0}^{1}x\, f(x) J_{0}(\lambda_kx) dx.\]
\end{theorem}
\section{Inverse source problem}
In this section, we consider the following inverse source problem:\\ 
 Find a pair of functions $\left\lbrace u(x,t), h(x) \right\rbrace $, in a  rectangular domain $\Omega,$ which satisfies the boundary value problem:
 \begin{equation}\label{sprob}
^c D^{\alpha}_{0|t}u=\dfrac{\partial^2 u}{\partial x^2}+\dfrac{1}{x}\dfrac{\partial u}{\partial x}+h(x),
\end{equation}

\begin{equation}
\underset{x\rightarrow 0}{\lim} \;x \;\dfrac{\partial u}{\partial x}=0,\quad u(1,t)=0,
\end{equation}

\begin{equation}\label{scond}
u(x,0)=g(x),\; u(x,T)=f(x).
\end{equation}
Here $g(x)$ and $f(x)$ are given functions and  $ \,^cD^{\alpha}_{0|t}$ stands for the Caputo fractional derivative of order $0<\alpha< 1$.
\subsection{Existence of solution }
Using the solution method described in the previous section, we seek solutions of the form 
\[ u(x,t)=\sum_{k=1}^{\infty} u_k(t) J_0(\lambda_k x), \quad h(x)=\sum_{k=1}^{\infty} h_k J_0(\lambda_k x), \]
where the coefficients $u_k(t)$ and $h_k$ are the unknowns to be determined. Substituting these representations into $(\ref{sprob})$ and $(\ref{scond})$, we obtain the following non-homogeneous fractional differential equation :
\begin{equation}\label{nfde}
 \,^c D^{\alpha}_{0|t}u_k(t)+\lambda_k^2  u_k(t)=h_k,
\end{equation}
with the boundary conditions :
\begin{equation}\label{nfdec}
u_k(0)=g_k,\quad u_k(T)=f_k,
\end{equation}
where $g_k$ and $f_k$ are the coefficients of Fourier-Bessel series for $g(x)$  and $f(x)$, respectively, which are given by
$$I_k=\dfrac{2}{J_{1}^{2}(\lambda_k)}\int_{0}^{1}x\, I(x) J_{0}(\lambda_k x) dx,\quad I= g,f.
$$
Solving $(\ref{nfde})$ we obtain (\cite{gm}, p. 243)
$$u_k(t)= C_k E_{\alpha}(-\lambda_k^2 t^\alpha)+\dfrac{h_k}{\lambda_k^2},$$
and using the conditions  $(\ref{nfdec})$, we deduce the following:
$$C_k=\dfrac{g_k-f_k}{1-E_{\alpha}(-\lambda_k^2 T^{\alpha})},\qquad h_k=\lambda_k^2(g_k-C_k).$$
Consequently, we have
\[\begin{array}{ll}
u(x,t)&=\displaystyle\sum_{k=1}^{\infty}\left( C_k E_{\alpha}(-\lambda_k^2 t^\alpha)+g_k-C_k\right)J_0(\lambda_k x), \\
&=g(x)-\displaystyle\sum_{k=1}^{\infty}\dfrac{{1-E_{\alpha}(-\lambda_k^2 t^{\alpha})}}{{1-E_{\alpha}(-\lambda_k^2 T^{\alpha})}}(g_k-f_k)J_0(\lambda_k x),\end{array}\]
and 
\[\begin{array}{ll}
h(x)&=\displaystyle \sum_{k=1}^{\infty}\lambda_k^2 (g_k-C_k)J_0(\lambda_k x),\\
&=g''(x)-\displaystyle\sum_{k=1}^{\infty}\dfrac{\lambda_k^2 (g_k-f_k)}{1-E_{\alpha}(-\lambda_k^2 T^{\alpha})}J_0(\lambda_k x).\end{array}\]
Using appropriate conditions on $f(x)$ and $g(x)$, as given in the theorem below, the convergence of the series expansions of $u(x,t), u_{x}(x,t), u_{xx}(x,t), ^c D^{\alpha}_{0|t}u(x,t)$ and $h(x)$ can be shown in a similar approach as discussed in the inverse initial problem in the previous section. 
\subsection{Main result}
The following theorem summarizes the main result of this section :
\begin{theorem}
Assume that the functions $g$ and $f$ are differentiable $4$ times such that 
\begin{itemize}
\item $g^{(i)}(0)= f^{(i)}(0)= 0,$ $(i=0,1,2,3)$;
\item $g^{(j)}(1)=f^{(j)}(1)=0,$ $(j=0,1,2)$;
\item $g^{(4)}(x)$ and $f^{(4)}(x)$ are bounded.
\end{itemize} 
 
Then, the inverse problem
\[\begin{array}{ll}
^c D^{\alpha}_{0|t}u=\dfrac{\partial^2 u}{\partial x^2}+\dfrac{1}{x}\dfrac{\partial u}{\partial x}+h(x),&0<x<1,\; 0<t<T,
\\
\underset{x\rightarrow 0}{\lim} \;x \;\dfrac{\partial u}{\partial x}=0,\quad u(1,t)=0,&0<t<T,\\
u(x,0)=g(x),\; u(x,T)=f(x),&0<x<1,
\end{array}\]
has a unique solution given by 
\[\begin{array}{rl}
u(x,t)=&g(x)-\displaystyle\sum_{k=1}^{\infty}\dfrac{1-E_{\alpha}(-\lambda_k^2 t^{\alpha})}{{1-E_{\alpha}(-\lambda_k^2 T^{\alpha})}}(g_k-f_k)J_0(\lambda_k x),\\
h(x)=&g''(x)-\displaystyle\sum_{k=1}^{\infty}\dfrac{\lambda_k^2 (g_k-f_k)}{1-E_{\alpha}(-\lambda_k^2 T^{\alpha})}J_0(\lambda_k x),
\end{array}\]
where $f_k$ and $g_k$ are given by
$$f_k=\dfrac{2}{J_{1}^{2}(\lambda_k)}\int_{0}^{1}x\, f(x) J_{0}(\lambda_k x) dx \quad \mbox{and} \quad g_k=\dfrac{2}{J_{1}^{2}(\lambda_k)}\int_{0}^{1}x\, g(x) J_{0}(\lambda_k x) dx.$$
\end{theorem}
\begin{description}
\item[Acknowledgements.] 
Authors acknowledge financial support from The Research Council (TRC), Oman. This work is funded by TRC under the research agreement no. ORG/SQU/CBS/13/030. Authors are also thankful to Dr. Erikin Karimov for useful discussions.
\end{description}
%%%%%%%%%%%%%%%%%%%%%%%%%%%%%%%%%%%%%%%%%%%%%%%%
%%%%%%%%%%%%%%%%%%%%%%%%%%%%%%%%%%%%%%%%%%%%%%%%
%%%%%%%%%%%%%%%%%%%%%%%%%%%%%%%%%%%%%%%%%%%%%%%%
%%%%%%%%%%%%%%%%%%%%%%%%%%%%%%%%%%%%%%%%%%%%%%%%


\begin{thebibliography}{99}
\bibitem{ams}T. Aleroev, M. Kirane and S. Malik, Determination of a source term for a time fractional diffusion equation with an integral type over-determining condition, Electronic Journal of Differential Equations, No. 270, pp. 1-16 (2013).

\bibitem{BB} M. Bertero, P. Boccacci,  Introduction to inverse problems in Imaging, IOP Publishing, Bristol, (1998).

\bibitem{furati} K. Furati, O.Iyiolaa and M. Kirane, An inverse problem for a generalized fractional diffusion, Journal of Applied Mathematics and Computation,
Vol. 249, p. 24-31 (2014).

\bibitem{gm} R. Gorenflo, F. Mainardi,  Fractional calculus: integral and differential equations of fractional order, Fractals and Fractional Calculus in continuum Mechanics : A. Carpinteri and F. Mainardi (eds), Springer Verlag, Wien and New York (1997),  pp 223-276.
\bibitem{hig}J.R. Higgins, Completeness and basis properties of sets of special functions, Cambridge University Press, 1977.
\bibitem{kirane and malik}M. Kirane , S. Malik, and M. Al-Gwaiz, An inverse source problem for a two dimensional time fractional diffusion equation with nonlocal boundary conditions, Journal of Mathematical Methods in Applied Sciences 36,  pp.1056-1069 (2013). 

\bibitem{initial} K. Masood, S. Messaoudi and F.D. Zaman, Initial inverse problem in heat equation with Bessel operator, International Journal of Heat and Mass Transfer 45, pp. 2959-2965 (2002).
\bibitem{moore} C.N. Moore, On the uniform convergence of the developments in Bessel functions, American Mathematical Society, vol 12, pp. 181-206 (1911).
\bibitem{Pet} T. St. Petrova, Application of Bessel's functions in the modelling of chemical engineering processes, Bulgarian Chemical Communications, vol 41, pp. 343-354 (2009).
\bibitem{Prabhakar} T. R. Prabhakar, A singular integral equation with a generalized Mittag-Leffler function in the kernal, Yokohama Math. J. 19, pp. 7-15 (1971).
\bibitem{Sch} M. G. Scherberg, The degree of convergence of a series of bessel functions, Trans. Am. Math. Soc. 35, pp. 172-183 (1933).
\bibitem{tolstov} G.P. Tolstov, Fourier series (translated by R.A. Silverman), Prentice Hall, Inc., Englewood Cliffs, N. J., (1962).

\bibitem{uhlmann1} Inside Out: Inverse Problems and Applications (MSRI Publications, Volume 47), Edited by G. Uhlmann, Cambridge University Press, Cambridge, (2003).

\bibitem{uhlmann2} Inverse Problems and Applications: Inside OutI II (MSRI Publications, Volume 60), Edited by G. Uhlmann, Cambridge University Press, Cambridge, (2013).

\bibitem{watson} G.N. Watson,  A treatise on the theory Of Bessel functions, Gambridge University Press, Second Edition  (1966).
 \bibitem{zhang} Y. Zhang and X. Xu, Inverse source problem for a fractional diffusion equation, Inverse Problems 27(3) (2011).
 
\end{thebibliography}
\end{document}